\documentclass[11pt]{article}
\usepackage{amsfonts}
\usepackage[mathscr]{eucal}
\usepackage{mathrsfs}
\usepackage{amsmath}
\usepackage{amssymb}
  \usepackage{paralist}
  \usepackage{graphics} 
  \usepackage{epsfig} 
\usepackage{graphicx}  \usepackage{epstopdf}
 \usepackage[colorlinks=true]{hyperref}
\hypersetup{urlcolor=blue, citecolor=red}

\newtheorem{theorem}{Theorem}[section]

\newtheorem{lemma}[theorem]{Lemma}

\newtheorem{definition}[theorem]{Definition}

\newcommand{\ep}{\varepsilon}

\newcommand{\De}{\Delta}
\newcommand{\La}{\Lambda}

\newcommand{\Om}{{\Omega}}

\newcommand{\al}{{\alpha}}

\newcommand{\si}{{\sigma}}

\newcommand{\om}{{\omega}}

\newcommand{\sL}{\mathscr{L}}

\newcommand{\R}{{\mathbb R}}

\newcommand{\cA}{\mathcal{A}}

\newcommand{\cC}{\mathcal{C}}
\newcommand{\cO}{\mathcal{O}}

\newcommand{\cE}{\mathcal{E}}
\newcommand{\cH}{\mathcal{H}}

\newcommand{\g}{{\nabla}}

\newcommand{\bu}{\bar{u}}
\newcommand{\pd}{\partial}

\newcommand{\arr}{\rightarrow}

\newcommand{\Lto}{\widehat L_2(\Om)}
\newcommand{\Hto}{\widehat H^2_{0}(\Om)}
\newcommand{\hch}{\widehat{\cH}}
\newcommand{\wH}{\widehat{H}}

\newcommand{\di}{{\rm div\, }}
\newcommand{\tv}{\tilde{v}}
\newcommand{\tw}{\tilde{w}}
\newcommand{\tu}{\tilde{u}}
\newcommand{\tU}{\widetilde{U}}
\newcommand{\tbu}{\tilde{\bar{u}}}
\newcommand{\wtE}{\widetilde{\cE}}

\begin{document}
\title{On an attractor for strong solutions to interactive fluid-plate system without rotational inertia.}

\author{Iryna Ryzhkova\thanks{e-mail: iryonok@gmail.com}\\
Humboldt Universit\"at zu Berlin}





\maketitle

\begin{abstract}
We study long-time dynamics of strong solutions to a 
 non-homogeneous coupled system consisting of linearized 3D
Navier--Stokes equations in a bounded domain  and the classical (nonlinear) full von Karman plate equations that account for both transversal and lateral displacements on a flexible   part  of the boundary. Rotational inertia of the filaments  of the plate is not taken into account.
Our main result is the existence of an attractor in the strong phase space provided lateral external loads are small enough.


\end{abstract}

\section{Introduction}
We study long-time behaviour of a coupled  system which describes 
interaction of a homogeneous viscous incompressible fluid
which  occupies a bounded domain $\cO$ and a thin (nonlinear) elastic plate.
The motion of the fluid is described
by the  linearized 3D Navier--Stokes equations and the  deformations of the plate
are described by the full von Karman plate model without rotational inertia.
Let us describe the problem in details.

Let $\cO\subset \R^3$  be a bounded domain  with a sufficiently smooth boundary $\partial\cO$. We assume that $\partial\cO=\Omega\cup S$, where
\begin{equation*}
\Om\subset\{ x=(x_1;x_2;0)\, :\,x'\equiv (x_1;x_2)\in\R^2\}
\end{equation*}
with a smooth contour $\Gamma=\partial\Om$, and $S$  is a  surface which lies in the halfspace $\R^3_- =\{ x_3\le 0\}$. The exterior normal to $\partial\cO$ is denoted  by $n$. Evidently, $n=(0;0;1)$ on $\Om$.

We consider the following {\em linear} Navier--Stokes equations in $\cO$ for the fluid velocity field $v=v(x,t)=(v^1(x,t);v^2(x,t);v^3(x,t))$ and the pressure $p(x,t)$:
\begin{align}
  & v_t-\nu\Delta v+\nabla p=G_{fl}\quad {\rm in\quad} \cO   \times(0,+\infty), \label{fl.1} \\
  &\di v=0 \quad {\rm in}\quad \cO  \times(0,+\infty), \label{fl.2}
\end{align}
where $\nu>0$ is the dynamical viscosity of the fluid and $G_{fl}$ is a volume force.
We supplement (\ref{fl.1}) and (\ref{fl.2}) with  the non-slip  boundary   conditions:
\begin{equation}\label{fl.4}
v=0 ~~ {\rm on}~S; \quad v\equiv(v^1;v^2;v^3)=(u^1_t;u^2_t;w_t) ~~{\rm on} ~ \Om,
\end{equation}
where $u=u(x,t)\equiv (u^1; u^2; w)(x,t)$ is the  displacement of the plate occupying $\Om$. Here $w$ stands for the transversal displacement, $\bu=(u^1;u^2)$ --- for the lateral (in-plane) displacements.

These boundary conditions describe influence of the plate on the fluid. Conversely,  the surface force $T_f(v)$ exerted by the fluid on the plate is equal to $Tn\vert_\Om$, where $n$ is the outer unit normal to $\pd\cO$ at $\Om$  and $T=\{T_{ij}\}_{i,j=1}^3$ is the stress tensor of the fluid,
\begin{equation*}
T_{ij}\equiv T_{ij}(v)=\nu\left(v^i_{x_j}+v^j_{x_i}\right)-p\delta_{ij}, \quad i,j=1,2,3.
\end{equation*}
Since  $n=(0;0;1)$ on $\Om$, we have that
\begin{equation}\label{surface-force}
    T_f(v)=(\nu(v^1_{x_3}+v^3_{x_1}); \nu(v^2_{x_3}+v^3_{x_2}); 2\nu \pd_{x_3}v^3 -p).
\end{equation}
\par
To describe the shell motion we use the full von Karman model which does not take into account rotational inertia of the filaments, but accounts for in-plane acceleration terms. We assume that
Young's modulus $E$ and Poisson's ratio $\mu\in(0,1/2)$ are such that $Eh=2(1+\mu)$, where $h$ is the thickness of the plate. The corresponding PDE system has the form
\begin{align}
	& w_{tt} + \De^2 w  =\di\left\{ \cC(P(u))\nabla w\right\}
	+ G_3- 2\nu \pd_{x_3}v^3 +p \label{pl_trans},\\
	& \bu_{tt}  = \di\left\{ \cC(P(u))\right\}  +
	     \left(
	       \begin{array}{c}
	        G_1-\nu (v^1_{x_3}+v^3_{x_1}) \\
	        G_2-\nu (v^2_{x_3}+v^3_{x_2})\\
	       \end{array}
	     \right),
	     \label{pl-plane}
\end{align}
where $G_{pl}=(G_1;G_2;G_3)$ is an external load and $\cC(P(u))$ is the stress tensor with
\begin{align*}
     & \cC(\epsilon)=2(1-\mu)^{-1}\left[\mu \, {\rm trace}\, \epsilon \cdot  I +(1-\mu)\epsilon \right], \qquad  &P(u)=\epsilon_0 (\bu)+  f(\nabla w), \\
     & \epsilon_0(\bu)=\frac 12 (\nabla \bu + \nabla ^T \bu), \qquad &f(s)=\frac 12 s \otimes s, \; s\in \R^2.
\end{align*}

This form of the full von Karman system was used earlier by many authors in the case when the   fluid velocity field $v$ is absent (see, e.~g.,~ \cite{KochLa_2002} and the references therein).

We impose the clamped boundary conditions on the plate
\begin{equation}
u^1|_{\pd\Om}=u^2|_{\pd\Om}=w|_{\pd\Om}=\left.\frac{\pd w}{\pd n'} \right|_{\pd\Om}=0, \label{pl-BC}
\end{equation}
where $n'$ is the outer normal to $\pd\Om$ in $\R^2$.

The resulting system \eqref{fl.1}-\eqref{pl-BC} is supplied with the initial data for the fluid velocity field $v=(v^1; v^2; v^3)$ and the plate displacement  $u=(u^1; u^2; w)$:
\begin{equation}
v\big|_{t=0}=v_0,~~ u\big|_{t=0}=u_0, ~~ u_t\big|_{t=0}=u_1.  \label{IC}
\end{equation}
Here $v_0=(v^1_0; v^2_0; v^3_0)$,  $u_j=(u^1_j; u^2_j; w_j)$, $j=0,1$, are given vector functions subjected to some compatibility conditions which we specify later.

We  note that  \eqref{fl.2} and \eqref{fl.4} imply the following
compatibility condition
\begin{equation}
\int_\Om w_t(x',t) dx'=0 \quad \mbox{for all}~~ t\ge 0.
\end{equation}
This condition fulfills when
\begin{equation}
    \int_\Om w(x',t) dx'=const \quad \mbox{for all}~~ t\ge 0 \label{ave_preserve}
\end{equation}
and can be interpreted as preservation of the volume of the fluid.

This fluid-structure interaction   model assumes that  large deflections of the elastic structure produce  small effect on the fluid. This  corresponds to the case when the fluid fills the container which is large in comparison with the size of the plate.

Mathematical studies of the problem
of fluid--structure interaction
in the case of viscous fluids and elastic plates/bodies
have a long history.
We refer to
\cite{AvaBuc2015,Chu_2010,ChuRyz2011,Ggob-jmfm08,Ggob-aa09,Ggob-mmas09}
and the references therein in the case of plates.

Long-time  dynamics for {\it nonlinear} plate-fluid
models was studied before in \cite{Chu_2010,ChuRyz2011,ChuRyzh2013-JDE, Ryzh_2018}.
The article \cite{Chu_2010} deals with
a class of fluid-plate interaction problems, when the plate, occupying $\Om$, oscillates in {\em longitudinal} directions only. A fluid-plate interaction model, accounting for purely transversal displacement of the plate without rotational inertia, was studied  in \cite{ChuRyz2011}. Existence and finite-dimensionality of an attractor were proved for both problems under the standard assumptions on non-linearities. 

Note,  that even in the linear case we cannot split
system \eqref{fl.1}--\eqref{IC} into two sets
of equations describing  longitudinal and
transversal plate movements separately, i.e.,
we cannot reduce the model under consideration to
the cases studied in \cite{Chu_2010,ChuRyz2011}.
For the detailed discussion one can see \cite{ChuRyzh2013-JDE}, Remark 1.1 (C). The paper mentioned above deals with a fluid-plate interaction model which describes simultaneous 
transversal and in-plane oscillations of the plate. In contrast to the model in the present paper, it accounts for rotational inertia of the plate filaments. To guarantee existence of an attractor for this model the authors had to assume mechanical dissipation in the transversal displacement equation. It is known, that the related {\em linear} system with rotational inertia accounted for lacks uniform stability (it is only strongly stable) \cite{AvaBuc2015}. That is why, probably, the model from \cite{ChuRyzh2013-JDE} in the absence of the mechanical dissipation does not possesses a global attractor.

Problem  \eqref{fl.1}--\eqref{IC} was first addressed in \cite{Ryzh_2018}.  In the system under consideration rotational inertia is neglected, thus $w_t$ has lower regularity ($L_2(\Om)$), then is the case of rotational inertia accounted for. Such regularity still allows us to prove existence of weak solutions satisfying  the energy inequality the same way as in \cite{ChuRyzh2013-JDE}, but  uniqueness of these solutions is still an open question. Therefore we had to resort to more smooth solutions. Well-posedness of  strong solutions to \eqref{fl.1}--\eqref{IC}  and their uniform stability in the case of zero external loads were proved in \cite{Ryzh_2018}. 
In this paper we proceed to the investigation of the attractor of the system with non-zero external loads. 

In Section \ref{sec:pre} we collect all the results about \eqref{fl.1}--\eqref{IC}  we need for further studies, and in Section \ref{sec:AB} we prove our main result on existence of an attractor. The main novelty of the paper is the proof of dissipativity in the strong phase space norm. The proof essentially relies on the obtained in \cite{Ryzh_2018} estimates of solution norms in the spaces of the form $L_2(0,T;H^s(\Om))$, where $H^s(\Om)$ is a space of smoother functions, than the phase space. We use finite difference estimates for Lyapunov function except of estimates for derivative because of that. Asymptotic smoothness is proved by the Ball's method.

\section{Preliminaries}\label{sec:pre}
In this section we introduce Sobolev type spaces we need,
provide some results concerning   the Stokes problem and collect the previous results on strong solutions to \eqref{fl.1}-\eqref{IC}, which we need to investigate its long-time behaviour.

\subsection{Spaces and notations}
To introduce Sobolev spaces we follow approach presented
in \cite{Triebel78}.
\par
Let $D$ be a sufficiently smooth domain  and $s\in\R$.
We denote  by $H^s(D)$ the Sobolev space of order $s$
on the set $D$ which we define as  a restriction (in the sense of distributions)
 of the
space $H^s(\R^d)$ (introduced via Fourier transform) to the domain $D$.
We define the norm in  $H^s(D)$ by the relation
\begin{equation*}
	\|u\|_{s,D}^2=\inf\left\{\|w\|_{s,\R^d}^2\, :\; w\in H^s(\R^d),~~ w=u ~~	\mbox{on}~~ D	\right\}.
\end{equation*}

We also use the notation $\|\cdot \|_{D}=\|\cdot \|_{0,D}$ and $(\cdot,\cdot)_D$
for the corresponding $L_2$ norm and  inner product. If the domain $D$ can be easily recognized from context we drop it in the notations.

We denote by $H^s_0(D)$ the closure of $C_0^\infty(D)$ in  $H^s(D)$
(with respect to  $\|\cdot \|_{s,D}$) and introduce the spaces
\[
H^s_*(D):=\left\{ f\big|_D\, :\; f\in H^s(\R^d),\;
{\rm supp }\, f\subset \overline{D}\right\},\quad s\in \R.
\]
Below we need them to describe
boundary traces on $\Om\subset\partial \cO$.
We endow the classes $H^s_*(D)$ with the induced norms
 $\|f \|^*_{s,D}= \| f \|_{s,\R^d}$
for $f\in H^s_*(D)$. It is clear that
\[
\|f \|_{s,D}\le \|f \|^*_{s,D}, ~~ f\in H^s_*(D).
\]
However, in general the norms $\|\cdot \|_{s,D}$ and
$\|\cdot \|^*_{s,D}$ are not equivalent.

Understanding adjoint spaces with respect to duality between
$C_0^\infty(D)$ and $[C_0^\infty(D)]'$
by Theorems 4.8.1 and 4.8.2 from \cite{Triebel78} we also have that
\begin{align*}
 [H^s_*(D)]'= H^{-s}(D),~ s\in\R, ~~~\mbox{and} ~~~
 [H^s(D)]' =  H_*^{-s} (D),~ s\in (-\infty,1/2).
\end{align*}
Below we also use the factor-spaces
$H^s(D)/\R$  with the naturally induced norm.
\medskip\par
To describe fluid velocity fields
we introduce the following spaces.
Let $\cC(\cO)$  be the class of $C^\infty$ vector-valued solenoidal (i.e., divergence-free) functions
on $\overline{\cO}$ which vanish in a neighborhood  of $S$ and  $\cC_0(\cO)$  be the class of $C^\infty_0(\cO)$ vector-valued solenoidal functions.
We denote by $X$ the closure of $\cC(\cO)$ with respect to the $L_2$-norm and
by $V$ the closure of $\cC(\cO)$ with respect to the $H^{1}(\cO)$-norm. Notations $V_0, \; X_0$ are used for the closure of $\cC_0(\cO)$ with respect to the $H^{1}(\cO)$-norm and $L_2(\cO)$-norm, respectively. One can see that
\begin{equation}\label{X-space}
X=\left\{ v=(v^1;v^2;v^3)\in [L_2(\cO)]^3\, :\; \di\, v=0,
\gamma_n v\equiv (v,n)=0~\mbox{on}~ S\right\};
\end{equation}
and
\begin{equation*}
    V=\left\{v=(v^1;v^2;v^3)\in [H^{1}(\cO)]^3\, :\; \di\, v=0,\; v=0~\mbox{on}~ S  \right\}.
\end{equation*}
We equip   $X$ and $X_0$ with the $L_2$-norm $\|\cdot\|_\cO$
and denote by $(\cdot,\cdot)_\cO$ the corresponding inner product.
We denote
\begin{equation*}
  V^\alpha=H^{1+\al}(\cO)\cap V, \qquad V^\alpha_0=H^{1+\al}(\cO)\cap V_0, \quad \al\ge 0.
\end{equation*}
The spaces $V^\al, \; V^\al_0$ are endowed  with the norm  $\|\cdot\|_{V^\al}= \|\nabla\cdot\|_{\al,\cO}$, and $V=V^0$.
For the details concerning  spaces of this type  we refer to \cite{temam-NS},
for instance.
\par
We also need the Sobolev spaces consisting of functions with zero average
on the domain $\Om$, namely
we consider the spaces
\begin{equation*}
	\widehat{L}_2(\Om)=\left\{u\in L_2(\Om): \int_\Om u(x') dx' =0 \right\}
\end{equation*}
and  $\widehat H^s(\Om)=H^s(\Om)\cap\widehat L_2(\Om)$ for $s>0$
with the standard $H^s(\Om)$-norm.
The notations   $\widehat H^s_*(\Om)$ and $\widehat H^s_0(\Om)$
have a similar meaning.

To describe plate displacement  we use the spaces
\begin{equation*}
W=H^1_0(\Om) \times H^1_0(\Om) \times H^2_0(\Om), \quad Y=L_2(\Om) \times L_2(\Om) \times \Lto
\end{equation*}
for weak solutions and
\begin{equation*}
W_s=H^1_0(\Om)\bigcap H^2(\Om)  \times H^1_0(\Om)\bigcap H^2(\Om) \times H^2_0(\Om)\bigcap  H^4(\Om), \quad Y_s=W
\end{equation*}
for strong solutions.
For weak solutions as a phase space we  use
\begin{equation*}
\cH= \left\{ (v_0;u_0;u_1)\in X\times W\times Y :\; v_0^3=u^3_1 ~\mbox{on}~ \Om\right\}, \qquad
\end{equation*}
with the standard product norm.

We also denote by $\hch$ a subspace in $\cH$ of the form
\begin{equation}\label{space-cH-hat}
\hch= \left\{ (v_0;u_0;u_1)\in \cH :\; w_0\in \wH^2_0(\Om)\right\},
\end{equation}
where $w_0$ is the third component of the  displacement
vector $u_0$. Phase space for strong solutions will be defined later.


\subsection{Stokes problem}

In further considerations we need some regularity  properties
of the terms responsible for fluid--plate interaction.
To this end we consider
 the following Stokes problem
\begin{align}
  -\nu\Delta v+\nabla p= g, \quad
   \di v=0 \quad {\rm in}\quad \cO; \nonumber
\\
 v=0 ~~ {\rm on}~S;
\quad
v=\psi=(\psi^1;\psi^2;\psi^3) ~~{\rm on} ~ \Om,\label{stokes}
\end{align}
where $g\in [L^2(\cO)]^3$ and $\psi\in [L^2(\Om)]^2 \times \Lto$ are given.
This type of boundary value problems for the Stokes equations
was studied by many authors  (see, e.g., \cite{temam-NS}
and references therein). We define operators we need and collect its properties in the following lemma.
\begin{lemma}\label{le:stokes}
	The following statements hold.
	\begin{itemize}
		\item [{ \bf (1)}] Let $g\in [H^{\sigma-1}(\cO)]^3$, and $\psi\in [H^{\sigma+1/2}_*(\Om)]^3$ with
			 $\int_\Om\psi^3(x')dx'=0$. Then for every $0\le \sigma\le  1$	problem \eqref{stokes} has a unique solution		$\{v;p\}$ in $[H^{1+\sigma}(\cO)]^3\times[ H^{\sigma}(\cO)/\R]$.
		\item[{\bf (2)}]
			We can define the operator $T_f: [H^{\sigma-1}(\cO)]^3\times [H^{\sigma+1/2}_*(\Om)]^3\arr [H^{\si-1/2}(\Om)]$ by the formula
			\begin{equation}\label{surface-force}
				T_f(g,\psi)=(\nu(v^1_{x_3}+v^3_{x_1}); \nu(v^2_{x_3}+v^3_{x_2}); 2\nu \pd_{x_3}v^3 -p),
			\end{equation}
			where $v$ is the solution to \eqref{stokes} with the right-hand side $g$ and boundary data $\psi$. This operator is linear and bounded between abovementioned spaces for every $0\le \sigma\le  1$.	
		\item [{ \bf (3)}] 
			We can define the linear operator $N_0 : [H^\si_*(\Om)]^2\times \widehat{H}^\si_*(\Om)\mapsto [H^{1/2+\si}(\cO)]^3\cap X$
			 by the formula
			\begin{equation}\label{fl.n0}
				N_0\psi=w ~~\mbox{iff}~~\left\{
					\begin{array}{l}
						 -\nu\Delta w+\nabla p=0, \quad  \di w=0 \quad {\rm in}\quad \cO;\\
						 w=0 ~~ {\rm on}~S;	\quad w=\psi ~~{\rm on} ~ \Om,
					\end{array}
				\right.
			\end{equation}
			This operator is linear and bounded between abovementioned spaces for every $0\le \sigma\le  1$.
	\end{itemize}
\end{lemma}
For the proof see \cite{ChuRyzh2013-JDE}.



\subsection{Well-Posedness Theorem}\label{sec:WP}
These definitions and theorems were formulated and proved in \cite{Ryzh_2018}. We remind them here for the reader's convenience.

To define weak (variational) solutions to \eqref{fl.1}--\eqref{IC}
we need the following class $\sL_T$ of test functions $\phi$ on $\cO$:
\begin{equation*}
\sL_T=\left\{\phi \left|\begin{array}{l}
\phi\in L_2(0,T; \left[H^1(\cO)\right]^3),\; \phi_t\in L_2(0,T;  [L_2(\cO)]^3),  \\
\di\phi=0,\; \phi|_S=0,\; \phi|_\Om=b=(b^1;b^2;d),  \\
d\in  L_2(0,T; \Hto),\; b^j \in L_2(0,T; H^1_0(\Om)), \; j=1,2,\; \\
d_t\in  L_2(0,T; \Lto), \; b^j_t \in L_2(0,T; L_2(\Om)), \; j=1,2.
\end{array}\right.\right\}.
\end{equation*}
We also denote $\sL_T^0=\{\phi\in \sL_T\, :\, \phi(T)=0\}$.
\begin{definition}[\cite{Ryzh_2018}]\label{de:solution}
{\rm
A pair of vector functions $(v(t);u(t))$ with $v=(v^1;v^2;v^3)$ and
$u=(u^1;u^2;w)$
is said to be  a weak solution to the problem \eqref{fl.1}--\eqref{IC}  on a time interval $[0,T]$ if
\begin{itemize}
    \item $v\in L_\infty(0,T;X)\bigcap L_2(0,T; V)$;
    \item $u \in L_\infty(0,T; H^1_0(\Om)\times H^1_0(\Om)\times H^2_0(\Om))$;
    \item $(u^1_t; u^2_t; w_t) \in L_\infty(0,T; L_2(\Om) \times L_2(\Om)\times \Lto)$  and   $u(0)=u_0$;
    \item $u_t \in L_2(0,T; \big[H^{1/2}_*(\Om)\big]^3)$ and the compatibility condition  $v(t)|_\Om=(u^1_t;u^2_t;w_t)(t)$ holds for almost all $t\in [0,T]$;
    \item for every $\phi\in \sL_T^0$  with $\phi|_\Om=b=(b^1;b^2;d)$ the following equality holds:
        \begin{equation}\label{weak_sol_def}
          \begin{split}
              -\!\int_0^T\!\!(v,\phi_t)_\cO dt & + \nu\!\int_0^T\!\!E(v,\phi)_{\cO}  dt -\!\int_0^T\!\! \big[(w_t,d_t)_\Om+(\bu_t,\bar{b}_t)_\Om\big]dt \\
               &  + \!\int_0^T\!\!(\De w, \De d)_\Om dt +  \!\int_0^T\!\!(\cC(P(u),\g \delta \otimes \g w + \epsilon_0(\bar{b}))_\Om dt = \\
               &(v_0, \phi(0))_{\cO} +(w_1,d(0))_\Om+ (\bu_1,\bar{b}(0))_\Om  + \\
               & \int_0^T(G_{fl}, \phi)_{\cO} dt +\int_0^T(G_{pl},b)_\Om dt
          \end{split}
        \end{equation}
        where  $\bu=(u^1;u^2)$, $\bar{b}=(b^1;b^2)$ and
        \begin{equation*}
            E(u,\phi)=\frac 12 \sum_{i,j=1}^3\left(v^j_{x_i} + v^i_{x_j} \right)\left(\phi^j_{x_i} + \phi^i_{x_j} \right).
        \end{equation*}
\end{itemize}
}
\end{definition}

The following theorem on existence of weak solutions can be proved the same way, as in \cite{ChuRyzh2013-JDE}.
\begin{theorem}[\cite{Ryzh_2018}] \label{th:wWP}
Assume that $U_0=(v_0;u_0;u_1)\in \cH$, $G_{fl}\in V'$,  $G_{pl}=(G_1; G_2; G_3)\in \left[H^{-1/2}(\Om)\right]^3$. Then  for any interval $[0,T]$ there exists a weak solution $(v(t); u(t))$ to \eqref{fl.1}--\eqref{IC}  with the initial data $U_0$. This solution possesses the following properties:
    \begin{equation*}\label{cont-ws}
    U(t;U_0)\equiv U(t)\equiv (v(t); u(t); u_t(t))\in L_\infty(0,T; X\times W \times Y),
    \end{equation*}
    \begin{equation}\label{smooth_est}
    ||w_t(t)||^2_{H^{1/2}_*(\Om)}+||u_t^1(t)||^2_{H^{1/2}_*(\Om)}+||u_t^2(t)||^2_{H^{1/2}_*(\Om)} \le C||\g v(t)||^2_\cO
    \end{equation}
    for almost all $t\in [0,T]$.

    The solution is bounded globally in $t$.
\end{theorem}

However, due to the strong supercriticality of the nonlinearity in the full von Karman equations we didn't manage to prove uniqueness of weak solutions. Therefore we build dynamical system in a space of more smooth functions.

To describe behaviour of the fluid component, we will use the space
\begin{equation}\label{ss-fl-space}
  X_s=V_0^1\oplus N_0(Y_s)
\end{equation}
where $N_0(Y_s)$ is the image of $Y_s$ under $N_0$ in $V^{1/2}$. That is, $v\in X_s$ means $v=\mathring{v}+N_0 u_1$, where $\mathring{v}\in V^1_0$ and $u_1\in Y_s$. The norm is defined as
\begin{equation*}
  ||v||^2_{X_s}=||\mathring{v}||^2_{V^1_0}+ ||N_0 u_1||_{3/2,\cO}.
\end{equation*}
For strong solutions we use the spaces
\begin{equation}\label{strong-spaces}
  \cH_s=\{(v_0, u_0, u_1)\in X_s\times W_s\times Y_s: \; v_0|_\Om=u_1 \}, \qquad \hch_s = \cH_s\cap \hch.
\end{equation}

\begin{definition}[\cite{Ryzh_2018}]
  A pair of vector functions $(v(t);u(t))$ is said to be a strong solution to problem \eqref{fl.1}--\eqref{IC}  on a time interval $[0,T]$ if it is a weak solution to this problem on $(0,T)$ and
  \begin{equation*}
    (v;u;u_t)(t)\in L_\infty(0,T; X_s\times W_s \times Y_s).
  \end{equation*}
\end{definition}

We also use the following  equivalent definition of variational solution. Taking in \eqref{weak_sol_def} $\phi(t)=\int_t^T\chi(\tau)d\tau\cdot \psi$, where $\chi$ is a smooth scalar function  and $\psi$ belongs to the space
\begin{equation*}\label{space-W}
\widetilde{V}=\left\{
\psi\in  V \left|  \;
\psi|_\Om=
\beta\equiv(\beta^1;\beta^2;\delta)\in H^1_0(\Om) \times H^1_0(\Om) \times \Hto \right. \right\}=V_0\oplus N(X_s),
\end{equation*}
one can see that the weak  solution $(v(t);u(t))$  satisfies the relation
\begin{equation}\label{weak_sol_d2}
\begin{split}
&\frac{d}{dt}[(v(t),\psi)_\cO + (w_t(t),\delta)_\Om+(\bu_t(t),\bar{\beta})_\Om]= \\
& -[ \nu E(v(t),\psi) +(\De w(t), \De \delta)_\Om + (\cC(P(u(t))),\g\delta \otimes \g w(t) +\epsilon_0(\bar{\beta}))_\Om] \\
& + (G_{fl}, \psi)_{\cO}  + (G_{pl},\beta)_\Om.
\end{split}
\end{equation}
for  all $t\in [0,T]$ and $\psi\in \widetilde{V}$  with $\psi\big\vert_\Om=\beta=(\beta^1;\beta^2;\delta)$ and $\bar{\beta}=(\beta^1;\beta^2)$.

Well-posedness of \eqref{fl.1}-\eqref{IC} in $\cH_s$ was proved in \cite{Ryzh_2018}.

\begin{theorem}[\cite{Ryzh_2018}] \label{th:sWP}
Assume that $U_0=(v_0;u_0;u_1)\in \cH_s$, $G_{fl}\in  X$,  $G_{pl}\in \left[H^{1/2}(\Om)\right]^3$. Then  for any interval $[0,T]$ there exists a unique strong  solution $(v(t); u(t))$ to \eqref{fl.1}--\eqref{IC}  with the initial data $U_0$. This solution possesses the following properties:
\begin{itemize}
\item it is continuous with respect to $t$ in the phase space, i.e.
    \begin{equation}\label{cont-ws:ss}
    U(t;U_0)\equiv U(t)\equiv (v(t); u(t); u_t(t))\in C(0,T; X_s\times W_s \times Y_s),
    \end{equation}
\item  there exists $C_T>0$, depending on $||U_0||_{\cH_s}$, such that for all $\ep>0$
\begin{equation} \label{smooth_est:ss}
    \begin{split}
         \int_0^T dt [ & ||w_t(t)||^2_{H^{5/2}(\Om)}  +
         ||\bu_t(t)||^2_{H^{3/2}(\Om)}+ ||\g v(t)||^2_{1-\ep,O}\\
          & + ||w(t)||^2_{H^{9/2}(\Om)}+||\bu(t)||^2_{H^{5/2}(\Om)}] \le C_{T}
    \end{split}
\end{equation}

\item The solutions depends continuously (in strong topology) on initial data in the space $\cH_s$.
\item The energy balance equality
\begin{equation}\label{lin_energy}
  \begin{split}
      \cE(v(t), u(t), u_t(t))+ & \nu \int_0^t E(v,v) d\tau  = \\
    & \cE(v_0, u_0, u_1) +\int_0^t(G_{fl},  v)_\cO d\tau +\int_0^t (G_{pl},u_t)_\Om d\tau
  \end{split}
\end{equation}
is valid for every $t>0$, where the energy functional $\cE$ is defined by
\begin{equation*}
\cE(v, u, u_t)=\frac12\left[\|v\|^2_\cO+ \|w_t\|^2_\Om+ \|\bu_t\|^2_\Om+\| \De w\|_\Om^2 +(\cC(P(u)),P(u))_\Om\right]
\end{equation*}
\end{itemize}
\end{theorem}


In the following lemma we collect additional properties of strong solutions to \eqref{fl.1}-\eqref{IC}, which we need to investigate its long-time behaviour.
\begin{lemma}\label{le:sWC}
	Let the conditions of Theorem \ref{th:sWP} hold. Then
	\begin{itemize}
		\item A strong solution depends continuously on initial data in weak topology of $\cH_s$.
		\item 'Energy equality of higher order' holds:
			\begin{equation}\label{diff-en-eq}		
				\wtE(U(t))+\int_{0}^{t} \nu E(\tv,\tv)d\tau=\wtE(U_0)   -  \frac{3}{2}\int_{0}^{t}(\cC(P(u,\tu)), \g\tw \otimes \g\tw)_\Om d\tau.
			\end{equation}
			where $\tv=v_t$, $\tbu=\bu_t$, $\tw=w_t$, $\tilde{u}=u_t$,
			\begin{multline}\label{diff-en}
			\wtE(U(t))= \frac{1}{2} \left[||\tv(t)||^2_\cO + ||\tu_t(t)||^2_\Om + ||\De\tu(t)||^2_\Om \right.\\
			\left. +(\cC(P(u,\tu)), P(u,\tu))(t) + 2(\cC(P(u)), \g\tw \otimes \g\tw)(t) \right],
			\end{multline}
			\begin{equation}\label{Puu_def}
			P(u(t),\tu(t))=\frac{d}{dt}P(u(t))=\ep_0(\tbu(t)) +\frac{1}{2} [\g w \otimes \g\tw + \g\tw\otimes \g w](t). 
			\end{equation}
		\item $v_t\in L^2(0,T;V)$ and
			\begin{equation}\label{vt_smooth}
				\int_0^T dt E(v_t,v_t) \le C(T, ||U_0||_{\cH_s}). 
			\end{equation}
		\item Strong solutions satisfy 
			\begin{multline}
				\int_{0}^{t}d\tau||\tw(\tau)||^2_{9/4,\Om} \le  2\delta \int_{0}^{t}d\tau E(\tv, \tv) + \\
				 \frac{C^w_1}{4\delta}\int_{0}^{t}d\tau \left[\frac{C^w_2}{4\gamma^2}\wtE(\tau)+||G|| + f(u)\right],  \label{w-d9/4-norm}
			\end{multline}	
			where $C_j^w$ are generic constants depending on  $\cH$-norm of initial data and $f(u)=f(||\De w||_\Om, ||\bu||_{1,\Om})$ is a generic function that behaves like $c_1||\De w||^\alpha_\Om + c_2||\bu||^\beta_{1,\Om}$ near zero with $\alpha,\beta>0$.
			$G$ means $(G_{fl}, G_{pl})$.
	\end{itemize}
\end{lemma}

The same way as in \cite[Theorem 3.3, Step 6]{ChuRyzh2013-JDE} we can prove that $S_t U_0^n \rightharpoonup S_t U_0$ in $\cH$ provided $U_0^n \rightharpoonup U_0$ in $\cH$. Now let $U_0^n \rightharpoonup U_0$ in $\cH_s$, therefore it bounded in $\cH_s$ and so is $S_t U_0^n$. Thus, it converges weakly to $S_t U_0$ in $\cH_s$ too.

Equality \eqref{diff-en-eq} was proved in \cite{Ryzh_2018} (see proof of Theorem 3.4, Step 1), and  \eqref{vt_smooth} is its direct consequence. The estimate \eqref{w-d9/4-norm} was proved in \cite{Ryzh_2018}.

\section{Asymptotic behaviour}\label{sec:AB}

In this section we study asymptotic behaviour of strong solutions to \eqref{fl.1}-\eqref{IC} within the dynamical system framework. Theorem \ref{th:sWP} and conservation of the average of the transversal displacement \eqref{ave_preserve} give us, that \eqref{fl.1}-\eqref{IC} generates dynamical systems $(\cH_s, S_t)$ and $(\hch_s, S_t)$. The first one cannot be dissipative because of \eqref{ave_preserve}, therefore we resort to study of the second one. Our main result is as follows.

\begin{theorem}
	Let conditions of Theorem \ref{th:sWP} hold true and in-plane external loads $(G_1,G_2)$ are small enough. Then the dynamical system $(\hch_s, S_t)$ possesses a compact global attractor.
\end{theorem}
In the proof we use the well-known scheme: dissipativity and asymptotic smoothness imply existence of an attractor \cite{Chu_2015}.

\subsection{Dissipativity}

For our purposes we extend the notion of dissipativity.

\begin{definition}
	Let $X\subset Y$ are Banach spaces and the embedding is continuous,  and $(X,S_t)$ be a dynamical system. We say that $(X,S_t)$ is $Y$-dissipative if there exists a set $B_0\subset X$ bounded in $Y$ such that for very set $D\subset X$ bounded in $Y$ there exists $T>0$ such that  for all $t>T$ $S_tD\subset B_0$.
\end{definition}

The following Lemma is the slightly modified result from \cite{Ryzh_2018}.
\begin{lemma}\label{le:est_below}
	Let $G\in [L_2(\Om)]^3$ and $||G_1||, ||G_2||$ are small enough. Then 
	\begin{equation*}
		||\De w||^2+(\cC(P(u)),P(u))_\Om +(G,u) \ge c_1||u||_W - C(G).
	\end{equation*}
	Let $||\g w||\le R$. Then
	\begin{equation}\label{est-u-quadr}
		||\De w||^2+(\cC(P(u)),P(u))_\Om \ge \frac{||u||_W^2}{C(R^2+2)}.
	\end{equation}
\end{lemma}

Using Lemma \ref{le:est_below} we can prove the following result in the way similar to Lemma 4.1 from \cite{Ryzh_2018}.
\begin{lemma}\label{le:weak_diss}
    Let $G\in [L_2(\Om)]^3$ and $||G_1||,  ||G_2||$ are small enough. Then  the set the stationary  points of problem \eqref{fl.1}-\eqref{IC} is bounded and the DS $(\hch_s, S_t)$ is $\hch$-dissipative. 
\end{lemma}    

Now we can prove dissipativity of $(\hch_s, S_t)$ in the strong norm.
\begin{lemma}\label{le:diss}
    Let $G_{fl}\in  X$,  $G_{pl}\in \left[H^{1/2}(\Om)\right]^3$ and let the dynamical system generated by \eqref{fl.1}--\eqref{IC} is $\hch$-dissipative with this right-hand side. Then the dynamical system is $\hch_s$-dissipative.
\end{lemma}

{\bf Proof}. We define Lyapunov function as
\begin{equation*}
    \La(U(t))= \wtE(U(t)) + \eta\left[ (\tu,\tu_t)_\Om(t) + (\tv, N_0\tu)_\cO(t) \right] + \bar{C},
\end{equation*}
where $\wtE$ is defined by \eqref{diff-en} and $\bar{C}$ is chosen to make $\La$ positive.  The constant $\eta>0$ will be chosen later.

Since the DS is $\hch$-dissipative, we can assume that initial data lie in the absorbing ball $B_a\subset\cH$ of the radius $R_d$. Thus, in this situation
\begin{equation}\label{equiv}
    c_1\wtE(U(t))-c_2 \le \La(U(t)) \le C_1\wtE(U(t)) +C_2.
\end{equation}
It was proved in \cite{Ryzh_2018}, that $\La$ is continuously differentiable with respect to $t$. The same way we obtain
\begin{multline}\label{Lyap_der}
  \frac{d}{dt} \La(U(t)) = -\nu E(\tv(t),\tv(t))+\eta||\tu_t(t)||^2_\Om + \Phi(U(t)) \\
  -\eta\left[ ||\De \tw(t)||^2_\Om + (C(P(u,\tu)),P(u,\tu))(t) \right],
\end{multline}
where
\begin{multline*}
    \Phi(U(t))=\eta\left[(\tv, N_0\tu_t)(t) - \nu E(\tv, N_0\tu)(t) - (C(P(u)),\g\tw\otimes\g\tw)(t) \right]- \\ \frac{3}{2}(C(P(u,\tu)),\g\tw\otimes\g\tw)(t).
\end{multline*}
In the estimates below we need to use estimate \eqref{w-d9/4-norm} for the $H^{9/4}(\Om)$-norm of $w$, which is integral with respect to $t$. Therefore we work with \eqref{Lyap_der} after integration with respect to time form $s$ to $t$.
Thus, \eqref{equiv} and \eqref{Lyap_der} imply that there exists $0<\omega=\omega(\eta)$ and $\bar{C}>0$ such that
\begin{multline}\label{Lyap_der1}
    \La(U(t))-\La(U(s)) +\om\int_{s}^{t} \La(U(\tau))d\tau +\frac{\nu}{2} \int_{s}^{t}E(\tv,\tv) d\tau \le \\
    \bar{C}(t-s) +\int_{s}^{t}\Phi(U(\tau))d\tau
\end{multline}

Let's estimate each term of $\int_{s}^{t}\Phi(U(\tau))d\tau$ now. First,
\begin{equation}\label{est2}
    \eta\left|\int_{s}^{t}(\tv, N_0\tu_t)_\cO d\tau\right| \le \eta C \int_{s}^{t}\left[E(\tv, \tv) + ||\tu_t||^2_\Om \right]d\tau \le
    \eta C \int_{s}^{t} E(\tv, \tv)d\tau
\end{equation}
due to trace theorem and Lemma~\ref{le:stokes}.
Second, using interpolation, energy equality \eqref{lin_energy} and the trace theorem  we obtain
\begin{multline}\label{est4}
    \eta \left|\int_{s}^{t} (C(P(u), \g\tw\otimes\g\tw))_\Om d\tau\right| \le C(R_d)\int_{s}^{t}||\tw||_\Om ||\De \tw||_\Om d\tau \le \\
    \eta C(R_d)\left[C(\delta) \int_{s}^{t} ||w_t||^2_\Om d\tau +\delta \int_{s}^{t} ||\De \tw||^2_\Om d\tau \right] \le\\
    \eta\left(\delta \int_{s}^{t} ||\De \tw||^2_\Om d\tau  + C(R_d,\delta)(t-s)\right).
\end{multline}
Then similarly to the previous estimate we have
\begin{multline}\label{est3}
  \eta\nu\left|\int_{s}^{t} E(\tv, N_0\tu) d\tau\right| \le  \eta\nu\int_{s}^{t} E(\tv,\tv)^{\frac 12}||\tu||^{\frac 14}_{1,\Om}||\tu||^{\frac 14}_{\Om} d\tau \le \\     
  \eta  \int_{s}^{t} E(\tv, \tv)d\tau + \eta\delta \int_{s}^{t} ||\tbu||^2_{1,\Om} d\tau +
  \eta\delta \int_{s}^{t} ||\tw||^2_{1,\Om}d\tau + \eta C(\delta,R_d)(t-s)
\end{multline}
It is left to estimate
\begin{equation*}
    \int_{s}^{t} (C(P(u,\tu)), \g\tw\otimes\g\tw)_\Om d\tau.
\end{equation*}
It consists of the terms of the form
\begin{equation*}
    B_3=\int_{s}^{t}\int_{\Om}\tbu^i_{x_j}\tw_{x_l}\tw_{x_k} dx d\tau, \quad B_4=\int_{s}^{t}\int_{\Om}w_{x_i}\tw_{x_j}\tw_{x_l}\tw_{x_k} dx d\tau.
\end{equation*}
Due to the boundary conditions and interpolation inequalities
\begin{multline*}
     |B_3| \equiv \left|\int_{s}^{t}\int_{\Om}\tbu^i\tw_{x_l x_j}\tw_{x_k} dx d\tau\right| \le \max_{\tau\in[s,t]}||\tbu^i(\tau)||_\Om\int_{s}^{t} ||\tw||_{2,\Om} ||\tw||_{1,\Om}d\tau \le \\
     C(R_d)\int_{s}^{t}||\tw||^{\frac 43}_{9/4,\Om} ||\tw||^{\frac 23}_\Om d\tau\le \epsilon \int_{s}^{t} ||\tw||^2_{9/4,\Om}d\tau + C(R_d, \epsilon)(t-s)
\end{multline*}
for every fixed $\epsilon>0$.

Let's estimate $|B_4|$. If $p,q>1$ and $1/p+1/q=1$, then
\begin{multline*}
  |B_4| \le
 C \int_{s}^{t}d\tau \left(\int_{\Om}|w_{x_i}|^p dx\right)^\frac 1p \left| \int_\Om |\tw_{x_k}|^q|\tw_{x_l}|^q|\tw_{x_j}|^q dx \right|^\frac 1q \le\\
  C\max_{\tau\in[s,t]}||w_{x_i}(\tau)||_{L^p(\Om)}\int_{s}^{t} d\tau ||\tw_{x_k}(\tau)||_{L^{3q}(\Om)} ||\tw_{x_l}(\tau)||_{L^{3q}(\Om)} ||\tw_{x_j}(\tau)||_{L^{3q}(\Om)}.
\end{multline*}
Thus, we have
\begin{equation*}
  |B_4|\le C_q \max_{\tau\in[s,t]}||\De w(\tau)||_\Om \cdot \int_{s}^{t} ||\g \tw||^3_{L^{3q}(\Om)}d\tau
\end{equation*}
for every $q>1$.
To estimate $||\g \tw||^3_{L^{3q}(\Om)}$ we use Gagliardo–Nirenberg interpolation inequality: in every smooth bounded domain $\Om\subset \R^2$
\begin{equation*}
    ||D^ju||_{L^p(\Om)} \le C ||D^m u||^\al_{L^r(\Om)}||u||^{1-\al}_{L^s{\Om}},
\end{equation*}
if
\begin{equation*}
    \frac 1p =\frac{j}{2} +\left(\frac 1r -\frac m2 \right)\al +\frac{1-\al}{s}, \quad \frac jm \le \al \le 1.
\end{equation*}
Choosing  $r=s=2, \; j=1,\; m=2$ we obtain
\begin{equation*}
    ||\g \tw||_{L^{3q}(\Om)} \le C ||D^2\tw||^\al||\tw||^{1-\al}, \quad \al=1-\frac{1}{3q} \ge 1/2 \mbox{ for } q>1.
\end{equation*}
Then, using interpolation between $H^{9/4}(\Om)$ and $L^2(\Om)$, we arrive to
\begin{equation*}
    ||\g \tw||^3_{L^{3q}(\Om)} \le C ||\tw||^{3-\frac 1q}_{2,\Om} ||\tw||^{\frac 1q}_\Om \le ||\tw||^\frac{24q-8}{9q}_{9/4,\Om} ||\tw||^\frac{8+3q}{9q}_\Om.
\end{equation*}
If $q<4/3$ then $\frac{24q-8}{9q}<2$ and we can use G\"older inequality and obtain
\begin{equation*}
    ||\g \tw||^3_{L^{3q}(\Om)} \le \ep||\tw||^2_{9/4,\Om} + C_\ep ||\tw||^r_\Om
\end{equation*}
for some $r$, therefore
\begin{equation*}
    |B_4| \le \ep \int_{s}^{t}||\tw||^2_{9/4,\Om}d\tau + C(\ep, R_d)(t-s).
\end{equation*}
Finally, using \eqref{w-d9/4-norm} with $\delta=2$, e.g, we obtain
\begin{multline*}
  \int_{s}^{t} (C(P(u,\tu),\g\tilde{w}\times\g\tilde{w})_\Omega d\tau \le \ep \int_{s}^{t}||\tw||^2_{9/4,\Om}d\tau + C(\ep, R_d)(t-s) \le  \\
  \ep \frac{\nu}{2} \int_{s}^{t} d\tau E(\tv,\tv) + \ep C(R_d)\int_{s}^{t}d\tau \wtE(U(\tau)) + C(\ep, R_d, G)(t-s).
\end{multline*}
Thus, we have the following estimate for $s<t$:
\begin{multline}
    \left|\int_{s}^{t} \Phi(\tau)d\tau\right| \le (\eta C+ \frac{\epsilon\nu}{2})\int_{s}^{t}  E(\tv, \tv)d\tau + \\
    (\eta\delta+\epsilon C(R_d))\int_{s}^{t} \wtE(U(\tau))d\tau + C(\epsilon, \eta, R_d, G)(t-s).
\end{multline}

First we chose $\eta$ and $\epsilon$ such that $\eta C+ \frac{\epsilon\nu}{2}<\frac{\nu}{2}$. We obtain  the certain value for $\om$.  Then we chose $\delta, \epsilon$ such that $(\eta\delta+\epsilon C(R_d))<\frac{\om}{2c_1}$, where $c_1$ is from \eqref{equiv}.
In the end we arrive
\begin{equation*}
    \La(U(t))-\La(U(s))+\frac\omega(\eta,\delta,\epsilon) 2\int_s^t    \La(U(\tau))d\tau \le C(\eta, \delta, \epsilon, R_d, G) (t-s).
\end{equation*}
We can divide the estimate by $t-s$ and pass to the limit when $t-s\arr+0$. This way we obtain the estimate for the left derivative for $\La$. Since $\La$ is continuously differentiable with respect to $t$, we have
\begin{equation*}
    \frac{d}{dt}\La(U(t)) +\frac{\omega}{2}\La(U(t))\le C(R_d, G),
\end{equation*}
which implies
\begin{equation*}
    \La(U(t))\le \La(0)e^{-\frac{\omega}{2}t} +C(R_d,G)
\end{equation*}
Lemma is proved.

\subsection{Asymptotic smoothness}
We  use the Ball's method to prove asymptotic smoothness  (see \cite{Ball2004} and also \cite{MRW}).
For  convenience we recall the abstract theorem (in a slightly relaxed form) from \cite{MRW}  which represents the main idea of the method.
\begin{theorem}[\cite{MRW}]\label{th:ball}
	Let  $S_t$ be a semigroup of  strongly
	continuous operators in some Hilbert space $F$.
	Assume that operators $S_t$ are also weakly continuous in $F$ and
	there exist  a number $\om>0$ and functionals  $\Lambda$, $L$ and $K$ on $F$
	such that the equality
	\begin{multline}
		\La(U(t)) +\int_s^t L(U(\tau))e^{-2\om(t-\tau)}d\tau= \\
		\La(U(s)) e^{-2\om(t-s)}+\int_s^t K(U(\tau))e^{-2\om(t-\tau)}d\tau \label{ball_eq}
	\end{multline}
	holds on the trajectories $U(t)$ of the dynamical system $(F,S_t)$.
	
	Let the	functionals possess the following properties:
	\begin{itemize}
		\item[{\bf (i)}] $\La:\, F \arr \R_+$ is a continuous bounded functional and if
		$\{U_j\}_j$ is bounded sequence in $F$ and  $t_j \arr +\infty$ is such that (a) $S_{t_j}U_j \rightharpoonup U$ weakly in $F$, and (b) $\limsup_{n\arr\infty}\La(S_{t_j}U_j) \le \La(U)$,
		then $S_{t_j}U_j \arr U$ strongly in $F$.
		
		\item[{\bf (ii)}]  $K:\; F \arr \R$ is 'asymptotically weakly continuous' in the sense that if  $\{U_j\}_j$ is bounded in $F$, and  $S_{t_j}U_j \rightharpoonup U$ weakly in $F$ as   $t_j \arr +\infty$, then $K(S_\tau U)\in L^{loc}_1(\R_+)$ and
		\begin{equation}\label{K-conv}
		\lim_{j\arr\infty}\int_0^t e^{-2\om(t-s)}K(S_{s+t_j}U_j) ds = \int_0^t e^{-2\om(t-s)}K(S_{s}U) ds, \quad \forall t>0.
		\end{equation}
		\item[{\bf (iii)}]  $L$ is 'asymptotically weakly lower semicontinuous' in the sense that if $\{U_j\}_j$ is bounded in $F$, $t_j \arr +\infty$, $S_{t_j}U_j \rightharpoonup U$ weakly in $F$, then $L(S_\tau U)\in L^{loc}_1(\R_+)$ and
		\begin{equation}\label{L-LowEst}
		\liminf_{j\arr\infty}\int_0^t e^{-2\om(t-s)}L(S_{s+t_j}U_j) ds \ge \int_0^t e^{-2\om(t-s)}L(S_{s}U) ds, \quad \forall t>0.
		\end{equation}
	\end{itemize}
	Then the dynamical system  $(F,S_t)$ is asymptotically smooth.
\end{theorem}

In what follows we denote
\begin{itemize}
	\item $U_j=(v_j, \bu_j^0, w_j^0, \bu_j^1, w_j^1)$ --- point in $\cH_s$ (initial data);
	\item $S_tU_j=U_j(t)=(v_j, \bu_j,w_j, \pd_t\bu_j, \pd_t w_j)(t)$ --- the strong solution to \eqref{fl.1}-\eqref{IC} with the initial conditions $U_j$;
	\item $\tilde{U}_j(t)=\pd_t U_j(t)=(\tv_j,\tbu_j, \tw_j, \pd_t\tbu_j, \pd_t \tw_j)(t)$ --- time derivative of the strong solution $U_j(t)$.
\end{itemize}
Since we will use 'energy relation' for time derivatives as the functional $\La$, we need to make a transition from values to it's derivatives and back. For convenience we formulate the first transition as a lemma.
\begin{lemma}\label{le:trans}
		\begin{enumerate}
			\item if $S_{t_j}U_j=U_j(t_j) \rightharpoonup U$ weakly in $\cH_s$ then $\tU_j(t_j) \rightharpoonup \tU$ weakly in $\cH$;
			\item if $S_{t_j}U_j=U_j(t_j) \arr U$ strongly in $\cH_s$ then $\tU_j(t_j) \arr \tU$ strongly in $\cH$;
		\end{enumerate}
\end{lemma}	
		We will prove the first statement. The second stetment can be proved the same way.
		let $S_{t_j}U_j \rightharpoonup  U=(v,\bu, w, \tbu, \tw)$ weakly in $\cH_s$. That is,
		\begin{align*}
			& \mathring{v}_j(t_j) \rightharpoonup  \mathring{v} \quad \mbox{ weakly in } V^1_0;\\
			&\bu_j(t_j) \rightharpoonup \bu \quad \mbox{ weakly in } [H^2(\Om)\bigcap H^1_0(\Om)]^2;\\
			& w_j(t_j) \rightharpoonup w \quad \mbox{ weakly in }  H^4(\Om)\bigcap H^2_0(\Om);\\
			& \pd_t\bu_j(t_j) \rightharpoonup \tbu \quad \mbox{ weakly in } [H^1_0(\Om)]^2;\\
			& \pd_tw_j(t_j) \rightharpoonup \tw \quad \mbox{ weakly in }  H^2_0(\Om);
		\end{align*}
		Prove the convergence of $\tU_j(t_j)$ in the corresponding space. Evidently,
		\begin{align*}
			& \tbu_j(t_j)= \pd_t\bu_j(t_j) \rightharpoonup \tbu \quad \mbox{ weakly in } [H^1_0(\Om)]^2;\\
			& \tw_j(t_j) = \pd_tw_j(t_j) \rightharpoonup \tw \quad \mbox{ weakly in }  H^2_0(\Om);\\
			& \tv_j(t_j)=\pd_t v_j(t_j)=\cA_s \mathring{v}_j(t_j)  \rightharpoonup  \cA_s\mathring{v} = \tv \quad \mbox{ weakly in } X_0;
		\end{align*}
		where $\cA_s$ is a Stokes operator. Further,
		\begin{equation*}
			\pd_t\tbu_j(t_j) =\pd_{tt}\bu_j(t_j)={\di}(\cC(P(u_j(t_j))))+
			\left(
				\begin{array}{c}
					G_1-\nu (\pd_{x_3}v_j^1+\pd_{x_1}v_j^3)(t_j) \\
					G_2-\nu (\pd_{x_3}v_j^2+\pd_{x_2}v_j^3)(t_j)\\
				\end{array}
			\right)
		\end{equation*}
		First, 
		\begin{equation*}
			P(u_j(t_j))=\frac 12 (\g\bu_j(t_j)+\g^T\bu_j(t_j)) + \frac 12 \g w_j(t_j) \otimes \g w_j(t_j).
		\end{equation*}
		The first term weakly converges in $[H^1_0(\Om)]^2$ to $\g\bu+\g^T\bu$. $\g w_j(t_j)$ converges strongly in $H^{3-\epsilon}$. Sinse  $H^{3-\epsilon}$ is a multiplicative algebra, the second term converges  strongly in $H^{3-\epsilon}$.
		Thus,
		\begin{equation*}
			{\di}(\cC(P(u_j(t_j)))) \rightharpoonup {\di}(\cC(P(u))) \quad \mbox{weakly in } [L_2(\Om)]^2.
		\end{equation*}
		Second, using Lemma \ref{le:stokes} p. (2) with $\sigma=1/2$,  we obtain that
		\begin{equation*}
		T_{f,in}(\tv,\tu)=
			\left(
				\begin{array}{c}
					-\nu (v^1_{x_3}+v^3_{x_1}) \\
					-\nu (v^2_{x_3}+v^3_{x_2})\\
				\end{array}
			\right)
		\end{equation*}
		 is a bounded (and thus continuous) linear operator from $X_0\times [H^1_0(\Om)]^3$ to $[L_2(\Om)]^2$. Thus, $T_{f,in}(\tv_j(t_j), \tu^j_t(t_j))$ converges weakly in $[L_2(\Om)]^2$.
		Consequently,
		\begin{equation*}
			\pd_t\tbu_j(t_j) =\pd_{tt}\bu_j(t_j) \rightharpoonup \hat{\bu} \quad \mbox{weakly in } [L_2(\Om)]^2.
		\end{equation*}
		Convergence of $\pd_t\tw_j(t_j)$ can be proved the same way.

Now we are ready to prove the main lemma of this subsection.
\begin{lemma}\label{le:as_smooth}
	Let $G_{fl}\in  X$,  $G_{pl}\in \left[H^{1/2}(\Om)\right]^3$. Then the dynamical  system $(\cH_s,S_t)$ is asymptotically smooth (see \cite[Definition 2.2.1]{Chu_2015}).
\end{lemma}
{\bf Proof}. Now we check validity of the hypotheses of Theorem~\ref{th:ball} in our case.
Strong continuity of $S_t$ in $\cH_s$ was proved in \cite[Theorem 3.4]{Ryzh_2018},  weak continuity of $S_t$ in $\cH_s$ was established in Lemma \ref{le:sWC}.

Now we introduce functionals from  Theorem \ref{th:ball}. The functional $\La$ is the Lyapunov function 
\begin{equation*}
\La(U(t))= \wtE(U(t)) + \eta\left[ (\tu,\tu_t)_\Om(t) + (\tv, N_0\tu)_\cO(t) \right] + \bar{C}.
\end{equation*}
Differentiating $\La$ with respect to $t$ \eqref{Lyap_der} and choosing $\om<\eta$ we arrive
\begin{multline}\label{ball_func}
	\frac{d}{dt} \La(U(t)) +2\om\La(U(t)) = \\
	-(\eta-\om)[||\De \tw(t)||^2 + (\cC(P(u,\tu)),P(u,\tu))(t)+ (\cC(P(u)),\g\tw \otimes \g\tw)(t)] + \\
	(\eta+\om)||\tu_t(t)||^2 +\om||\tv(t)||^2 +\eta[(v,N_0\tu_t)(t)-\nu E(\tv,N_0 \tu)(t)] -\nu E(\tv, \tv)(t) -\\
	\frac 32(\cC(P(u,\tu)),\g\tw \otimes \g\tw)(t) + 2\om\eta[(\tu,\tu_t)(t) + (\tv, N_0\tu)(t)].
\end{multline}
We denote
\begin{multline}\label{LDef}
	L(U(t))=\\
	(\eta-\om)[||\De \tw(t)||^2 + (\cC(P(u,\tu)),P(u,\tu))(t)+ (\cC(P(u)),\g\tw \otimes \g\tw)(t)] + \\
	-(\eta+\om)||\tu_t(t)||^2-\om||\tv(t)||^2 +\nu E(\tv, \tv)(t)-\frac 32(\cC(P(u,\tu)),\g\tw \otimes \g\tw)(t),
\end{multline}
\begin{multline}\label{KDef}
	K(U(t))=\eta[(v,N_0\tu_t)(t)-\nu E(\tv,N_0 \tu)(t)] -\\
	\frac 32(\cC(P(u,\tu)),\g\tw \otimes \g\tw)(t) + 2\om\eta[(\tu,\tu_t)(t) + (\tv, N_0\tu)(t)].
\end{multline}
Thus, multiplying \eqref{ball_func} by $e^{-2\om(t-\tau)}$ and integrating with respect to $\tau$ we obtain \eqref{ball_eq}.

Let's now check the properties of the functionals. 

Property (i) of $\La$. Continuity on $\cH_s$ is evident, because $\La$ can be represented as a sum of the norm and the compact part. Boundedness follows from Lemma \ref{le:diss}. Let $t_j\arr +\infty$ and $\{U_j\}_{j=1}^\infty$ is a bounded in $\cH_s$ sequence such that $S_{t_j}U_j=U_j(t_j) \rightharpoonup U$ weakly in $\cH_s$ and $\limsup_{n\arr\infty} \La(U_j(t_j)) \le \La(U)$. We need to prove that $S_{t_j}U_j=U_j(t_j) \arr U$ strongly in $\cH_s$. 

First we prove that  
\begin{equation}\label{ner}
	\liminf_{n\arr\infty} \La(U_j(t_j)) \ge \La(U).
\end{equation}
From Lemma \ref{le:trans}  we have
\begin{equation*}
	\pd_t U_j(t_j) \rightharpoonup \widetilde{U} \quad  \mbox{ weakly in }  \cH.
\end{equation*}
Thus, $\tu_j(t_j) \arr \tu $ strongly in $H^{1-\ep}(\Om)$ and therefore 
\begin{equation}\label{conv1}
	\left[ (\tu_j,\pd_t\tu_j)(t_j) + (\tv_j, N_0\tu_j)(t_j) \right]  \arr \left[ (\tu,\tu_t) + (\tv, N_0\tu) \right].
\end{equation}
For the terms in $\wtE$, which are norms, the estimate like \eqref{ner} is evidently take place. 
Now we consider the term $(\cC(P(u_j)),\g\tw_j \otimes \g\tw_j)(t_j)$. Since $\g w(t_j)\arr w$ strongly in $H^{3-\ep}$ and  $H^{3-\ep}$ is a multiplicative algebra, $\cC(P(u_j))(t_j)\arr \cC(P(u))$  strongly in $H^{1-\ep}$. 
$\g \tw_j \rightharpoonup \g \tw$ weakly in $H^1(\Om)$, thus strongly in $H^{1-\ep}\subset L^4(\Om)$. Therefore 
\begin{equation}\label{conv2}
	(\g \tw_j \otimes \g \tw_j)(t_j) \arr \g \tw \otimes \g \tw \qquad \mbox{ strongly in }	L_2(\Om)
\end{equation}
and 
\begin{equation}\label{conv3}
	(\cC(P(u_j)),\g\tw_j \otimes \g\tw_j)(t_j) \arr (\cC(P(u)),\g \tw \otimes \g \tw).
\end{equation}
It is left to consider the term $(\cC(P(u_i,\tu_j)),P(u_i,\tu_j))_\Om$. Recollecting definition \eqref{Puu_def}, we have
\begin{equation}\label{conv4}
	\liminf_{j\arr \infty} \frac{1}{4} (\cC(\ep_0(\tbu_j)),\ep_0(\tbu_j))(t_j)  \ge \frac{1}{4} (\cC(\ep_0(\tbu)),\ep_0(\tbu)), 
\end{equation}
since this term is an equivalent norm in $[H^1_0(\Om)]^2$. From \eqref{conv2} it follows that
\begin{gather}
	((\cC(\ep_0(\tbu_j)), \g w_j\otimes \g\tw_j)(t_j) \arr 	((\cC(\ep_0(\tbu)), \g w\otimes \g\tw), \label{conv5}\\
	(\cC(\g w_j\otimes \g\tw_j), \g w_j\otimes \g\tw_j)(t_j) \arr (\cC(\g w\otimes \g\tw), \g w\otimes \g\tw). \label{conv6}
\end{gather}

Thus, we have proved \eqref{ner} and therefore $\La(U_j(t_j)) \arr \La(U)$. This implies that  $||S_{t_j}\tU_j||_\cH \arr ||\tU||_\cH$ and thus the convergence $S_{t_j}\tU_j \arr \tU$ is strong in $\cH$. 

Provided with this convergence we have the following consequences of \eqref{fl.1}-\eqref{pl-plane}. Equation \eqref{fl.1} implies $\mathring{v}_j(t_j)=\cA_s^{-1} \tv_j(t_j) \arr \cA_s^{-1} \tv$ strongly in $H^2(\cO)\bigcap X_0$, and, since $N_0 \tu_j(t_j) \arr N_0 \tu$ strongly in $N_0(W)$,
\begin{equation*}
	v_j(t_j)= \mathring{v}_j(t_j) +N_0 \tu_j(t_j)  \arr \cA_s^{-1} \tv +  N_0 \tu=v \quad \mbox{ strongly in } X_s.
\end{equation*}
Equations \eqref{pl_trans}, \eqref{pl-plane} together with Proposition \ref{le:stokes} give that
\begin{equation*}
	u_j(t_j) \arr u \qquad \mbox{strongly in } W_s.
\end{equation*}

Let us now check property (ii) of functional $K$. Let again $t_j\arr +\infty$ and $\{U_j\}$ is a bounded in $\cH_s$ sequence such that $S_{t_j}U_j=U_j(t_j) \rightharpoonup U$ weakly in $\cH_s$. Then Lemma \ref{le:trans} implies 
\begin{equation*}
	\pd_t U_j(t_j) \rightharpoonup \widetilde{U} \quad  \mbox{ weakly in }  \cH.
\end{equation*} 
Using \eqref{conv1} and \eqref{conv2} we can easily verify that
\begin{multline*}
	\eta(v_j,N_0\pd_t\tu_j)(t_j) -\frac 32(\cC(P(u_j,\tu_j)),\g\tw_j \otimes \g\tw_j)(t_j) +\\
	 2\om\eta[(\tu_j,\pd_t\tu_j)(t_j)+ (\tv_j, N_0\tu_j)(t_j) ] \arr \\
	\eta(v,N_0\tu_t) -\frac 32(\cC(P(u,\tu)),\g\tw \otimes \g\tw) + \\
	2\om\eta[(\tu,\tu_t) + (\tv, N_0\tu) ].
\end{multline*}
It is left to investigate the term  $E(\tv,N_0 \tu)$. First, Lemma \ref{le:sWC} implies $\tv_j(s+t_j)\in L_2(0,T;V)$ as functions of $s$ and they are bounded in this space by $C(T,R)$. Since $\tv_j(\cdot+t_j) \rightharpoonup \tv(\cdot)$ in $L_2(0,T;X)$ due to Lemma \ref{le:trans} we get
\begin{equation}\label{conv7}
	\tv_j(\cdot+t_j)\arr \tv(\cdot) \quad \mbox{ weakly in }L_2(0,T;V)
\end{equation}
Finally, $N_0 \tu_j(s+t_j) \rightharpoonup N_0 \tu(s)$ weakly in $H^{3/2}\bigcap X$ and thus strongly in $V$ for every fixed $s$. The property (ii) of Theorem \ref{th:ball} is proved.

Let us proceed to the property (iii) of functional $L$. We prove estimate \eqref{L-LowEst} for each term in $L$ separately.  \eqref{L-LowEst} is true for the first term $||\De \tw||^2_\Om$ since it is a norm. For the second term $(\cC(P(u,\tu)),P(u,\tu))_\Om$ \eqref{L-LowEst} follows from \eqref{conv4}-\eqref{conv6} and \eqref{conv3} implies \eqref{L-LowEst} for the third term $(\cC(P(u)),\g\tw \otimes \g\tw)_\Om$ of $L$. For the last term $(\cC(P(u,\tu)),\g\tw \otimes \g\tw)_\Om$ estimate \eqref{L-LowEst} follows from \eqref{conv2} and weak convergence of $\cC(P(u_j(t_j),\tu_j(t_j)))$ in $L_2(\Om)$. The sixth term $ E(\tv, \tv)$ satisfies \eqref{L-LowEst} because of weak convergence property \eqref{conv7}. The rest terms $-||\tu_t||^2_\Om$ and $ -\om||\tv||^2_\cO$ have wrong signs, so we need to prove strong convergence for them. Let us investigate $\tu_{tt}$ and $\tv_t$.
Differentiating \eqref{weak_sol_d2} with respect to $t$ we obtain that $\tu$ and $\tv$ satisfy
\begin{equation}\label{approx_sol_diff}
\begin{split}
&[(\tv_t(t),\psi)_\cO + (\tw_{tt}(t),\delta)_\Om+(\tbu_{tt}(t),\bar{\beta})_\Om]= \\
& -[ \nu E(\tv(t),\psi) +(\De \tw(t), \De \delta)_\Om + (\cC(P(u(t),\tu(t))),\g\delta \otimes \g w(t) +\ep_0(\bar{\beta}))_\Om] \\
& + (\cC(P(u(t))),\g\delta\otimes\g\tw(t)))_\Om,
\end{split}
\end{equation}
for every $\psi\in V_0\oplus N_0(Y_s)$. Thus, 
\begin{equation*}
	\pd_t\tv_j (\cdot +t_j) \rightharpoonup \tv_t(\cdot) \quad \mbox{ in } L_2(0,T; (V_0\oplus N_0(Y_s))').
\end{equation*}
Since $V_0\oplus N_0(Y_s) \Subset X\cap H^{1-\ep}(\cO) \subset X\subset (V_0\oplus N_0(Y_s))'$, the Aubin's lemma gives that
\begin{equation*}
	\tv_j (\cdot +t_j) \arr \tv_t(\cdot) \quad \mbox{ in } L_2(0,T; X\cap H^{1-\ep}(\cO)), \quad \ep\in(0,1].
\end{equation*}
Boundary conditions \eqref{fl.4} and the trace theorem imply
\begin{equation*}
	\pd_{t}\tu_j(\cdot +t_j) \arr \tw_{tt}(\cdot) \quad \mbox{ in } L_2(0,T; L^2(\Om)).
\end{equation*}
The Lemma is proved.

\section*{Acknowlegements}

The Author is grateful to Ukrainian armed forces for possibility to continue her research work and to Germany for hospitality. This work was financially supported by Volkswagen foundation (projects "From Modelling and Analysis to Approximation" and "Dynamic Phenomena in Elasticity Problems").

\end{document}